\documentclass[12pt]{article}
\usepackage{mathtools}
\usepackage[utf8]{inputenc}
\usepackage{amssymb,amsmath,amsfonts,eucal,mathrsfs,amsthm} 
\usepackage{xcolor}
\newtheorem{theorem}{Theorem}
\newtheorem{proposition}[theorem]{Proposition}
\newtheorem{lemma}[theorem]{Lemma}

\newtheorem{remark}[theorem]{Remark}

\theoremstyle{definition}
\newcommand{\R}{\mathbb{R}}
\newcommand{\Z}{\mathbb{Z}}
\newcommand{\Q}{\mathbb{Q}}
\newcommand{\Sf}{\mathbb{S}}
\newcommand{\C}{\mathbb{C}}
\newcommand{\Hy}{\mathbb{H}}

\newcommand{\spa}{\mbox{span}}

\newcommand{\Ric}{\mbox{Ric}}

\newcommand{\trace}{\mbox{tr\,}}

\def\<{{\langle}}
\def\>{{\rangle}}

\def\CP{\mathord{\mathbb C}\mathord{\mathbb P}}

\def\n{\nabla}

\def\a{\alpha}

\def\be{\begin{equation} }
\def\ee{\end{equation} }

\newcommand\blfootnote[1]{\begingroup
\renewcommand\thefootnote{}\footnote{#1}
\addtocounter{footnote}{-1}
\endgroup}
\begin{document}

\title{Ricci pinched compact submanifolds\\ in spheres}
\author{Marcos Dajczer and Theodoros Vlachos}
\date{}
\maketitle

\begin{abstract}
We investigate the topology of the compact submanifolds in 
round spheres that satisfy a lower bound on the Ricci 
curvature depending only on the length of the mean curvature 
vector of the immersion. Just in special cases, the limited 
strength of the assumption allows some strong additional 
information on the extrinsic geometry of the submanifold. 
\end{abstract}
\blfootnote{\textup{2020} \textit{Mathematics Subject 
Classification}: 53C20, 53C40, 53C42, 55R25.}

This paper is devoted to analyzing the topology and some aspects 
of the geometry of compact submanifolds in round spheres. Our 
starting point is the seminal paper by Lawson and Simons \cite{LS} 
which considered isometric immersions $f\colon M^n\to\Sf^{n+m}$ 
of compact $n$-dimensional manifolds into the unit sphere with 
codimension $m\geq 1$. They proved that under certain upper
bounds, given in terms of the second fundamental form of the 
submanifold, specific integral homology groups must vanish. 

The strategy in \cite{LS} goes as follows: From the work of 
Federer and Fleming it is known that stable minimal currents 
represent homology classes with integer coefficients. Then 
the results are achieved by ruling out the existence of 
these currents for certain dimensions. Since in a homology 
class the area can be minimized, this trivializes the 
corresponding integral homology groups for these 
dimensions.

In its strongest version, the bound required by Lawson and 
Simons is quite technical and  has no clear geometric meaning; 
see Theorem \ref{ls} below. However, they also showed that it 
can be replaced by an upper bound on the norm of the second 
fundamental form of $f$, which in their words, amounts for the 
submanifold to be ``reasonably untwisted". In fact, they prove 
that if $p,q$ are positive integers with $p+q=n$ such that 
$\|\a\|^2<\min\{pq,2\sqrt{pq}\}$ then $H_p(M^n,\Z)=H_q(M^n,\Z)=0$.
Moreover, the result is sharp since for all $p,q$ there exist 
embeddings of $\Sf^p\times\Sf^q$ in $\Sf^{n+1}$ such that 
$\|\a\|^2=2\sqrt{pq}$.

It was shown by Vlachos \cite{V1,V2} that the limitations  
of the strongest bound considered by Lawson and Simons 
can be overcome by assuming instead a strict lower bound for 
the Ricci curvature at any point of $M^n$. This bound is 
given in terms of an integer $2\leq k\leq n/2$ and the 
length of the mean curvature vector of the submanifold. 
Moreover, when $k=n/2$ if $n$ is even or $k=(n-1)/2$ if $n$ 
is odd, the conclusion is that $M^n$ is a homology sphere. 
Notice that for the class of minimal submanifolds the bound 
condition is purely intrinsic. 
\vspace{1ex}

The purpose of this paper is to investigate the topology and,
to some degree, the extrinsic geometry of the submanifolds when 
the requirement on the Ricci curvature in Vlachos' results is 
no longer limited to be strict. Our answer to this natural 
question shows that the class of submanifolds satisfying the 
less restricted bound is quite richer and the proofs 
are substantially more involved.  
Nevertheless, it is important to emphasize that several arguments 
presented in this paper significantly leverage insights advanced 
by Vlachos. In addition, we show that the aforementioned results 
by Vlachos still hold if the bound is required to be strict just 
at some point. Besides, it also follows from this paper that some 
results in \cite{ER} can be improved to now sharp statements. 
\vspace{1ex}

For an isometric immersion $f\colon M^n\to\Sf^{n+m}$
with second fundamental form with  values in the normal 
bundle $\a_f\colon TM\times TM\to N_fM$ its  mean 
curvature vector field is given by $\mathcal H=(1/n)\trace\a_f$.
Then the length of $\mathcal H$ is denoted by $H$.

We recall that  a vector in the normal bundle 
$\eta\in N_fM(x)$  at $x\in M^n$ is named a 
\emph{Dupin principal normal} of $f\colon M^n\to\Sf^{n+m}$ 
if the tangent vector subspace
$$
E_\eta(x)=\left\{X\in T_xM\colon\alpha_f(X,Y)
=\<X,Y\>\eta\;\,\text{for all}\;\,Y\in T_xM\right\}
$$
is at least two dimensional. That dimension is called the 
\emph{multiplicity} of $\eta$.

\begin{theorem}\label{thm1} Let $f\colon M^n\to\Sf^{n+m}$, 
$n\geq 4$, be an isometric immersion of a compact manifold. 
For an integer $2\leq k\leq n/2$ assume that at any point of 
$M^n$ we have that
\be
\Ric_M\geq b(n,k,H)=
\frac{n(k-1)}{k}+\frac{n(k-1)H}{2k^2}\big(nH+\sqrt{n^2H^2
+4k(n-k)}\big)\tag{$\ast$}
\ee 
where $\Ric_M$ stands for the (not normalized) Ricci curvature 
of $M^n$. Then $M^n$ is simply connected 
and hence orientable, and one of the following holds:\vspace{1ex} 

\noindent $(i)$ The homology groups satisfy that
$$
H_i(M^n;\Z)=H_{n-i}(M^n;\Z)=0\;\,\text{for all}\;\,1\leq i\leq k,
$$  
and if $k<n/2$ then $H_{n-k-1}(M^n;\Z)$ is torsion free. 
This is necessarily the case if $(*)$ is strict at some point. 
\vspace{1ex}

\noindent $(ii)$ The homology groups satisfy that 
$$
H_i(M^n;\Z)=H_{n-i}(M^n;\Z)=0\;\,\text{for all}\;\,1\leq i\leq k-1,
$$
$H_k(M^n;\Z)\neq0$ and $H_{n-k}(M^n;\Z)$ is torsion free. 
If for $k=2$ we assume in addition that $H\neq 0$ everywhere, 
then there is a Dupin principal normal vector 
$\eta\in N_fM(x)$ at any $x\in M^n$ collinear with 
$\mathcal H(x)$ when $H(x)\neq 0$. Moreover, its multiplicity 
$\ell(x)$ at $x\in M^n$
satisfies $k\leq\ell(x)\leq n-k-1$ if $k<n/2$ and $\ell(x)=k$ if $n$ 
is even and $k=n/2$. Furthermore, we have that
$$
\|\eta\|=\lambda(n,k,H)
=\frac{1}{2k}\big(nH+\sqrt{n^2H^2+4k(n-k)}\big)
$$ 
and that $\Ric_M(X)=b(n,k,H)$ for any unit vector $X\in E_\eta(x)$. 
\end{theorem}

The example given below by part $(i)$ of Proposition \ref{proj} 
for $m>2$ has the homology as in part $(ii)$ of the above result.
Moreover, the same example shows that an additional 
assumption that has been made for $k=2$ to conclude 
the existence of a Dupin principal normal vector 
is necessary.\vspace{1ex}

A submanifold $f\colon M^n\to\Sf^{n+m}$ is said to carry  a 
\emph{Dupin principal normal vector field} $\eta\in\Gamma(N_fM)$ 
if $\eta$ is nowhere vanishing and $\eta(x)$ is at any point  
$x\in M^n$ a Dupin principal normal with constant multiplicity 
$\ell\geq 2$. For any $f$ as in part $(ii)$ in the above result 
there is an open dense subset of $M^n$ such that along each connected 
component the submanifold carries a Dupin principal 
normal vector field.\vspace{1ex}

Throughout the paper that $(*)$ is \emph{satisfied
with  equality} at $x\in M^n$ means that there exists a unit 
vector $X\in T_xM$ such that $\Ric_M(X)=b(n,k,H)$. If 
otherwise, we call the inequality $(*)$ \emph{strict} 
at $x\in M^n$.
Notice that any compact submanifold that has strict 
inequality $(*)$ at any point continues trivially to 
have it after a sufficiently small smooth deformation.
\vspace{1ex}

The pinching condition $(*)$ was inspired by the 
generalized Clifford torus, that is, the standard embedding  
of $\mathbb{T}^n_p(r)=\Sf^p(r)\times\Sf^{n-p}(\sqrt{1-r^2})$,
$2\leq p\leq n-2$, into the unit sphere $\Sf^{n+1}$, where 
$\Sf^p(r)$ denotes the sphere of dimension $p$ and radius 
$r<1$. Computations in \cite{V2} yield that $(*)$ holds if 
$k=p$ and $(k-1)/(n-2)\leq r^2\leq k/n$,
and then  with equality.

In particular, the torus $\mathbb{T}^{2m+1}_m(r)$ with 
$r=\sqrt{m/(2m+1)}$ is a minimal hypersurface in $\Sf^{2m+2}$ 
that satisfies $(*)$ for $k\leq m$ with equality attained  
only for $k=m$ and for unit vectors tangent to the first 
factor. Moreover, it holds that $H_m(\Sf^m\times\Sf^{m+1};\Z)
=H_{m+1}(\Sf^m\times\Sf^{m+1};\Z)=\Z$. Hence the topological 
description as well as the characteristics of the Dupin 
principal normal given by part $(ii)$ of Theorem \ref{thm1} 
are both optimal. 

The next two results analyse the special cases when the 
value of $k$ in the inequality $(*)$ is the largest allowed 
in Theorem \ref{thm1}, namely, when $k=n/2$ if $n$ is 
even or $k=(n-1)/2$ if $n$ is odd. 
Unexpectedly, it turns out that between the two cases
there are striking differences in arguments and results.
Recall that a submanifold is called substantial if it is 
not contained in a proper totally geodesic submanifold of 
the ambient space.
Thus the assumption in the next result that the submanifold 
is substantial is just a technical requirement.

\begin{theorem}\label{thm2} Let $f\colon M^n\to\Sf^{n+m}$, 
$n\geq 4$, be a substantial isometric immersion of a compact 
manifold of even dimension. Assume that
\be\label{strict1}
\Ric_M\geq(n-2)\big(1+H^2+H\sqrt{1+H^2}\big)
\ee
at any point of $M^n$. Then one of the following holds:
\vspace{1ex}

\noindent $(i)$ $M^n$ is homeomorphic to $\Sf^n$ and this 
is necessarily the case if at some point of $M^n$
the inequality \eqref{strict1} is strict.\vspace{1ex}

\noindent $(ii)$ The submanifold is the minimal embedding 
of the torus $\mathbb{T}^n_{n/2}(1/\sqrt{2})$ into 
$\Sf^{n+1}$ given above.
\vspace{1ex}

\noindent $(iii)$ The submanifold is the minimal embedding 
$\psi\colon\CP^2_{4/3}\to\Sf^7$ given by part $(i)$
of Proposition \ref{proj}.
\end{theorem}

The topological aspect in the following result is the best 
possible from the homology point of view.

\begin{theorem}\label{thm3} Let $f\colon M^n\to\Sf^{n+m}$,
$n\geq 5$, be an isometric immersion of a compact manifold
of odd dimension.  Assume that
\be\label{ineq2}
\Ric_M\geq\frac{n(n-3)}{n-1}\Big(1+\frac{H}{n-1}
\big(nH+\sqrt{n^2H^2+n^2-1}\big)\Big)
\ee
is satisfied at any point of $M^n$. Then  we have the following:
\vspace{1ex}

\noindent $(i)$ $M^n$ is homeomorphic to $\Sf^n$ and 
this is necessarily the case if at some point of $M^n$ 
the inequality \eqref{ineq2} is strict. 
\vspace{1ex}

\noindent $(ii)$ If for $n=5$ we assume further that $H\neq 0$ 
at any point, then $M^n$ is diffeomorphic to the total space
$\mathsf{E}$ of a sphere bundle
$\Sf^k\hookrightarrow\mathsf{E}\xrightarrow{p}L^{k+1}$
with the base homeomorphic to $\Sf^{k+1}$. Moreover, either
\begin{itemize}

\item[(a)] If $n=5\,\text{or}\,13$ then $M^n$ is homeomorphic to 
$\Sf^k\times\Sf^{k+1}$. The homology of $M^n$ for the other 
dimensions is isomorphic to
the one of  $\Sf^k\times\Sf^{k+1}$, this being
necessarily the case if $n=4r+1$, or
\item[(b)] The homology of $M^n$ is $H_k(M^n,\Z)=\Z_q$
for some $q>1$, satisfies  $H_0(M^n,\Z)=H_n(M^n,\Z)=\Z$ 
and is trivial in all other cases.
\end{itemize}
\end{theorem}

The part $(i)$ of the above results relate with Theorem 4.1 
by Xu and Gu \cite{GUXU} where under the strict inequality 
$\Ric_M>(n-2)(1+H^2)$ it was established that the submanifold 
is homeomorphic to a sphere. In the minimal case, a complete 
classification for $\Ric_M\geq n-2$ is due to Ejiri \cite{E}
and is used as part of the proof of our Theorem \ref{thm2}.
\vspace{1ex}

Theorem \ref{thm3} reduces the topological classification of
the submanifolds of odd dimension satisfying the Ricci pinching 
condition to the classification of the total spaces of sphere 
bundles over spheres. By the bundle classification theorem 
given in \S $18.5$ of \cite{St} the equivalence classes of the 
$\Sf^k$ bundles over $\Sf^\ell$ are in one-to-one correspondence 
with the homotopy group $\pi_{\ell-1}(SO(k+1))$. Various sphere 
bundles over spheres that are not product bundles, but have the 
homotopy groups and homology structure of the products, are 
exhibited in $\S$26 of \cite{St}. It is worth mentioning that there 
is just one nontrivial equivalence class of $\Sf^4$ bundle 
over $\Sf^5$ since $\pi_4(SO(5))\cong\Z_2$. 
Similar situation holds for the $\Sf^{10}$ bundle over $\Sf^{11}$, 
$\Sf^{12}$ bundle over $\Sf^{13}$ and $\Sf^{14}$ bundle over 
$\Sf^{15}$. The $3$-sphere bundles over $\Sf^4$ are in one-to-one 
correspondence with $\pi_3(SO(4))\cong\Z\oplus\Z$.  
Among the total spaces of them there is the exotic 
$7$-sphere of Milnor as part of a family of $15$ examples. 
\vspace{1ex}

We are indebted to Wolfgang Ziller for his helpful comments.

\section{The pinching condition}

The result by Lawson and Simons \cite{LS} aforementioned
was strengthened by Elworthy and Rosenberg \cite[p.\ 71]{ER} 
by not requiring the bound to be strict at all points of 
the submanifold. In this section, we first state their result
and then we analyse the relation between the condition required 
in their theorems and our pinching assumption.

\begin{theorem}{\em(\cite{ER},\cite{LS})}\label{ls} 
Let $f\colon M^n\to\Sf^{n+m}$, $n\geq 4$, 
be an isometric immersion of a compact manifold and  $p$ an integer
such that $1\leq p\leq n-1$. Assume that at any point $x\in M^n$ 
and for \emph{any} orthonormal basis $\{e_j\}_{1\leq j\leq n}$ 
of $T_xM$ the second fundamental form 
$\alpha_f\colon TM\times TM\to N_fM$  satisfies
\be
\Theta_p=\sum_{i=1}^p\sum_{j=p+1}^n\big(2\|\a_f(e_i,e_j)\|^2
-\<\a_f(e_i,e_i),\a_f(e_j,e_j)\>\big)\leq p(n-p) \tag{$\#$}.
\ee
If at a point of $M^n$ the inequality $(\#)$ is strict for
any orthonormal basis, then there 
are no stable $p$-currents and the homology groups satisfy 
$H_p(M^n;\mathbb{Z})=H_{n-p}(M^n;\mathbb{Z})=0$.
\end{theorem}

The following is the main ingredient in the proofs of 
our theorems.

\begin{proposition}\label{prop}
Let $f\colon M^n\to\Sf^{n+m}$ be an isometric immersion 
satisfying the inequality $(*)$  at $x\in M^n$ for an integer 
$k$ such that $2\leq k\leq n/2$. Then at $x$ the following 
assertions hold:
\vspace{1ex}

\noindent $(i)$ The inequality $(\#)$ is satisfied for any 
$1\leq p\leq k$ and is strict for $p<k$ and any orthonormal 
basis $\{e_j\}_{1\leq j\leq n}$ of $T_xM$.   
Moreover, if the inequality $(*)$ is strict then 
also is $(\#)$ for $p=k$ and any orthonormal basis 
$\{e_j\}_{1\leq j\leq n}$ of $T_xM$. 
\vspace{1ex}

\noindent $(ii)$  If $k=2$ assume further that $H(x)\neq 0$. 
If equality holds in $(\#)$ for an orthonormal basis 
$\{e_j\}_{1\leq j\leq n}$ of $T_xM$ then $p=k$ and there is 
a Dupin principal normal vector $\eta$ at $x$ with 
$\|\eta\|=\lambda(n,k,H)$, that is collinear with ${\cal H}(x)$ if 
$H(x)\neq 0$, such that $\spa\{e_1,\dots,e_k\}\subset E_\eta(x)$.

Conversely, 
if there is a Dupin principal normal $\eta$ at $x$ of 
multiplicity at least $k$ with $\|\eta\|=\lambda(n,k,H)$
and is collinear with ${\cal H}(x)$ if 
$H(x)\neq 0$, then equality holds in $(\#)$ for $p=k$ 
and a certain orthonormal basis $\{e_j\}_{1\leq j\leq n}$
for which $\spa\{e_1,\dots,e_k\}\subset E_\eta(x)$.

Moreover, in the above situation $\Ric_M(X)=b(n,k,H)$ 
for any unit vector $X\in\spa\{e_1,\dots,e_k\}$. 
\vspace{1ex}

\noindent $(iii)$ If $H(x)=0$ and $k=2$, then equality 
holds in $(\#)$ for an orthonormal basis 
$\{e_j\}_{1\leq j\leq n}$ of $T_xM$ if and only if $(a)$ $p=2$, 
$(b)$ there is $\eta\in N_fM(x)$ with $\|\eta\|\leq\lambda(n,2,0)$  
satisfying
$$
\alpha_f(X,Y)=\<X,Y\>\eta
$$
for any $X,Y\in\spa\{e_1,e_2\}$ and $(c)$ $\Ric_M(X)=b(n,2,0)$ 
for any unit vector $X\in\spa\{e_1,e_2\}$. 
\end{proposition}

To prove the above result we need the following lemma.

\begin{lemma}\label{lemma}
Let $f\colon M^n\to\Sf^{n+m}$ be an isometric immersion 
satisfying the inequality $(*)$ at $x\in M^n$ for an integer 
$k$ such that $2\leq k\leq n/2$. Then the following assertions 
hold:\vspace{1ex}

\noindent $(i)$ If $H(x)\neq 0$ then the shape operator $A_1$ 
associated to $\xi_1=\mathcal{H}/H$ satisfies 
\be\label{ineq}
\mu(n,k,H)\leq\<A_1X,X\>\leq\lambda(n,k,H)
\ee
for any unit vector $X\in T_xM$, where 
\begin{align}
\lambda(n,k,H)\label{lambda}
&=\frac{1}{2k}\big(nH+\sqrt{n^2H^2+4k(n-k)}\big),\\
\mu(n,k,H)\label{mu}
&=\frac{1}{2k}\big(n(2k-1)H-\sqrt{n^2H^2+4k(n-k)}\big)
\end{align}
and $\<A_1X,Y\>=\<\a_f(X,Y),\xi_1\>$ for any $X,Y\in T_xM$.

\noindent $(ii)$ If $H(x)=0$ then $A_\xi$ associated to any 
unit vector $\xi\in N_fM(x)$ satisfies 
$$
\mu(n,k,0)\leq\<A_\xi X,X\>\leq\lambda(n,k,0)
$$
for any unit vector $X\in T_xM$. 
\end{lemma}

\proof $(i)$ First we recall that the Gauss equation implies that 
the Ricci curvature of a submanifold $f\colon M^n\to\Sf^{n+m}$ 
in the direction of a unit tangent vector $X\in T_xM$ at 
$x\in M^n$ is given by 
\be\label{ric}
\Ric_M(X)=n-1+\sum_{\a=1}^m(\mathrm{tr}A_\a)
\<A_\a X,X\>-\sum_{\a=1}^m\|A_{\a}X\|^2, 
\ee
where $A_{\a}, 1\leq\a\leq m$, are the corresponding shape 
operators of $f$ associated to an orthonormal basis 
$\{\xi_\a\}_{1\leq\a\leq m}$ of the normal space 
$N_fM(x)$ of $f$ at $x$.\vspace{1ex}  

Assume first that $H(x)\neq 0$. Lemma $2$ in \cite{V2} 
gives that the  mean curvature vector satisfies 
$n^2(k-2)H^2<4(n-k)$ with the use of which we obtain 
that $\mu(n,k,H)<\lambda(n,k,H)$.
\vspace{1ex}  

Let $\{\xi_\a\}_{1\leq\a\leq m}$ be an orthonormal basis 
of the normal space of $f$ at $x$ such that 
$\mathcal{H}=H\xi_1$ and denote by 
$A_{\alpha },1\leq\alpha\leq m$, the corresponding shape 
operators. Notice that $\trace A_1=nH$ and $\trace A_\a=0$ 
if $\a\geq 2$. From \eqref{ric} we have for any 
unit vector $X\in T_xM$ that
$$
\Ric_M(X)=n-1+nH\<A_1X,X\>-\|A_1X\|^2
-\sum_{\alpha\geq 2}\|A_\alpha X\|^2.
$$
Using $(*)$ we obtain that
$$
\|A_1 X\|^2-nH\<A_1X,X\> +b(n,k,H)-n+1\leq0
$$
and thus
$$
\<A_1X,X\>^2-nH\<A_1X,X\> +b(n,k,H)-n+1\leq0.
$$
Since \eqref{lambda} and \eqref{mu} are the roots of 
the  polynomial
\be\label{poli}
P(t)=t^2-nHt+b(n,k,H)-n+1
\ee
then \eqref{ineq}  has been proved. 
\vspace{1ex}

\noindent $(ii)$ Now assume that $H(x)=0$ and let $\xi$ be an arbitrary 
unit normal vector of $f$ at $x\in M^n$. 
Let $\{\xi_\a\}_{1\leq\a\leq m}$ be an orthonormal 
normal basis at $x$ with $\xi_1=\xi$ and corresponding 
shape operators $A_\a,1\leq\a\leq m$. 
Then for any unit vector 
$X\in T_xM$ we have  from \eqref{ric} that
$\Ric_M(X)=n-1-\sum_\a\|A_\a X\|^2$.
Then $(*)$ gives $\|A_1X\|^2+b(n,k,0)-n+1\leq 0$
and thus $\<A_1X,X\>^2\leq(n-k)/k$. Hence
$$
\mu(n,k,0)\leq \<A_\xi X,X\>\leq\lambda(n,k,0),
$$
which completes the proof.\vspace{2ex}\qed

\noindent\emph{Proof of Proposition \ref{prop}:}
$(i)$ Let $\{\xi_\a\}_{1\leq\a\leq m}$ be an orthonormal basis 
of $N_fM(x)$ were we agree that if $H(x)\neq 0$ then
$\mathcal{H}(x)=H(x)\xi_1$. Let $A_{\alpha}$, 
$1\leq\a\leq m$, be the corresponding shape operators 
and $\{e_j\}_{1\leq j\leq n}$ an orthonormal basis of 
$T_xM$. For simplicity, we denote 
$\a_{ij}=\alpha_f(e_i,e_j)$, $1\leq i,j\leq n$. Then
\begin{align}\label{a}
\Theta_p&=\;2\sum_{i=1}^{p}\sum_{j=p+1}^n\|\a_{ij}\|^2-n
\sum_{i=1}^p\<\a_{ii},\mathcal{H}\>
+\sum_{i,j=1}^p\<\a_{ii},\a_{jj}\>\nonumber\\
=&\;2\sum_{i=1}^p\sum_{j=p+1}^n\sum_\a\<A_{\alpha }e_i,e_j\>^2-nH
\sum_{i=1}^p\<A_1e_i,e_i\>
+\sum_\a\big(\sum_{i=1}^p\<A_\alpha e_i,e_i\>\big)^2\nonumber\\
\leq&\; 2\sum_{i=1}^p\sum_{j=p+1}^n\sum_\a\<A_{\alpha }e_i,e_j\>^2-nH
\sum_{i=1}^p\<A_1e_i,e_i\>
+p\sum_\a\sum_{i=1}^p\<A_{\alpha}e_i,e_i\>^2, 
\end{align}
where for the inequality we used that
\be\label{a-}
\big(\sum_{i=1}^p\<A_\alpha e_i,e_i\> 
\big)^2\leq p\sum_{i=1}^p\<A_{\alpha}e_i,e_i\> ^2. 
\ee

We first assume that $p=1$. Then \eqref{a} implies that
$$
\sum_{j=2}^n\big(2\|\a_{1j}\|^2-\<\a_{11},\a_{jj}\>\big)
\leq 2\sum_\a\|A_{\alpha}e_1\|^2
-nH\<A_1e_1,e_1\>.
$$
Then appealing to \eqref{ric} we obtain
$$
\sum_{j=2}^n\big(2\|\a_{1j}\|^2-\<\a_{11},\a_{jj}\>\big)
\leq 2(n-1-\Ric_M(e_1))+nH\<A_1e_1,e_1\>.
$$
From  $(*)$ it follows that
\be\label{ric1}
\Ric_M(e_i)\geq \frac{n(k-1)}{k}(1+H\lambda(n,k,H))
\;\,\text{for all}\;\,1\leq i\leq n.
\ee
On the other hand, Lemma \ref{lemma} yields
\be\label{A}
\<A_1e_i,e_i\>\leq\lambda(n,k,H)
\;\,\text{for all}\;\,1\leq i\leq p.
\ee
Then, from \eqref{ric1}, \eqref{A} and since $k\geq 2$ 
we obtain
$$
\sum_{j=2}^n\big(2\|\a_{1j}\|^2-\<\a_{11},\a_{jj}\>\big)
\leq\frac{2(n-k)}{k}
+\frac{2-k}{k}nH\lambda(n,k,H)<n-1.
$$

Next we assume that $p\geq 2$. We have 
\begin{align}\label{b-}
2\sum_{i=1}^p&\sum_{j=p+1}^n
\sum_\a
\<A_{\alpha}e_i,e_j\>^2
+p\sum_{i=1}^p\sum_\a
\<A_{\alpha}e_i,e_i\>^2\nonumber\\
&\leq p\sum_{i=1}^p\sum_{j=p+1}^n
\sum_\a\<A_{\alpha }e_i,e_j\>^{2}+p\sum_{i=1}^p\sum_\a
\<A_{\alpha}e_i,e_i\>^2\nonumber\\
&\leq p\sum_{i=1}^p\sum_\a\|A_{\alpha}e_i\|^2. 
\end{align}
Thus \eqref{a} implies that
$$
\Theta_p\leq p\sum_{i=1}^p\sum_\a\|A_{\alpha}e_i\|^2
-nH\sum_{i=1}^p\<A_1e_i,e_i\>.
$$
On the other hand, using \eqref{ric}, 
\eqref{ric1} and \eqref{A} it follows that
\begin{align}\label{c}
p\sum_{i=1}^p\sum_\a&\|A_{\a}e_i\|^2
-nH\sum_{i=1}^p\<A_1e_i,e_i\>\nonumber\\
&= p\sum_{i=1}^p(n-1-\Ric_M(e_i))
+(p-1)nH\sum_{i=1}^p\<A_1e_i,e_i\>\nonumber\\
&\leq\frac{p^{2}(n-k)}{k}
+\frac{p(p-k)}{k}nH\lambda(n,k,H). 
\end{align}
Since $p\leq k$ we have that 
\be\label{e}
\Theta_p\leq\frac{p^{2}(n-k)}{k}\leq p(n-p)
\ee
and $(\#)$ has been proved.

We already saw that that the inequality $(\#)$ is strictly
for $p=1$. If $p\geq 2$ and  $p<k$ the same conclusion 
follows from \eqref{e}. Finally, if $p=k$ and the 
inequality $(*)$ is strict then also \eqref{e} becomes 
strict, and this completes the proof of part $(i)$.
\vspace{1ex}

\noindent $(ii)$ Assume that at $x\in M^n$ equality holds 
in $(\#)$  for an orthonormal basis $\{e_j\}_{1\leq j\leq n}$ 
of $T_xM$.  Then part $(i)$ yields $p=k$ and 
\be\label{ricb}
\Ric_M(e_i)=b(n,k,H)\;\,\text{for all}\;\,1\leq i\leq k.
\ee
Moreover, all inequalities from \eqref{a-} 
to \eqref{e} become equalities. In particular, we have 
from \eqref{a-} that there exists $\rho_\a$ independent of $i$
such that
\be\label{r}
\<A_\a e_i,e_i\>=\rho_\a\;\,\text{for all}
\;1\leq i\leq k,\;\,1\leq\a\leq m.
\ee
From \eqref{b-} it follows that
\be\label{k}
(k-2)\<A_\a e_i,e_j\>=0\;\,\text{for all}
\;\,1\leq i\leq k,\;k+1\leq j\leq n,\;1\leq\a\leq m, 
\ee
and
\be\label{ii}
\<A_\a e_i,e_{i'}\>=0\;\,\text{for all}\;\,1\leq i\neq i' 
\leq k,\;1\leq\a\leq m,
\ee
whereas \eqref{A} and \eqref{c} give
\be\label{la}
H(\<A_1e_i,e_i\>-\lambda(n,k,H))=0
\;\,\text{for all}\;\,1\leq i\leq k.
\ee

\noindent\emph{Case I}. Assume that $H(x)\neq 0$. 
Then \eqref{la} implies that 
\be\label{n1}
\<A_1e_i,e_i\>=\lambda(n,k,H)\;\,\text{for all}\;\,1\leq i\leq k.
\ee
Since $\lambda(n,k,H)$ is a root of the polynomial given by 
\eqref{poli} it follows from \eqref{n1} that
$$
\<A_1e_i,e_i\>^2-nH\<A_1e_i,e_i\>+b(n,k,H)-n+1=0.
$$
Then \eqref{ric} and \eqref{ricb} yield 
$$
\sum_{\a\geq 2}\|A_\a e_i\|^2=\<A_1e_i,e_i\>^2-\|A_1e_i\|^2\leq 0.
$$
Hence
$$
A_\alpha e_i=0\;\,\text{for all}\;\,2\leq\alpha\leq m,\;1\leq i\leq k,
$$
and 
$$
A_1e_i=\lambda(n,k,H)e_i\;\,\text{for all}\;\,1\leq i\leq k.
$$
Therefore 
$$
\alpha_f(X,Y)=\<X,Y\>\eta\;\,\text{where}\;\,\eta=\lambda(n,k,H)\xi_1
$$
for any $X\in\spa\{e_1,\dots,e_k\}$ and $Y\in T_xM$.
\vspace{1ex}

\noindent\emph{Case II}. Assume that $H(x)=0$ in 
which case \eqref{la} is trivial. Since $k\geq 3$ by assumption, 
it then follows from \eqref{r}, \eqref{k} and \eqref{ii} that 
\be\label{above}
A_\a e_i=\rho_\a e_i\;\,\text{for all}\;\,
1\leq i\leq k,\;1\leq\a\leq m.
\ee
Hence $\eta=\sum_\a\rho_\a\xi_\a$ is a Dupin principal normal
satisfying $\spa\{e_1,\dots,e_k\}\subset E_\eta(x)$.
Then from \eqref{lambda}, \eqref{ric}, \eqref{ricb} and 
\eqref{above} we obtain 
\begin{align*}
\|\eta\|^2&=\sum_\a\rho_\a^2=\sum_\a\|A_\a e_i\|^2
=n-1-\Ric_M(e_i)=n-1-b(n,k,0)\\
&=(n-k)/k=\lambda^2(n,k,0)
\end{align*}
as wished.
\vspace{1ex}

Conversely, let $\{e_j\}_{1\leq j\leq n}$ be an orthonormal 
basis of $T_xM$  such that 
$\spa\{e_1,\dots,e_k\}\subset E_\eta(x)$.
Then let $\{\xi_\a\}_{1\leq\a\leq m}$ be a orthonormal basis 
of $N_fM(x)$ such that $\eta=\lambda(n,k,H)\xi_1$. Hence
$$
A_1e_i=\lambda(n,k,H)e_i\;\;\text{and}\;\;
A_\a e_i=0\;\,\text{for any}\;\,1\leq i\leq k,\;2\leq\a\leq m.
$$
Using that $\eta$ and $\mathcal H(x)$ are colinear if 
$H(x)\neq 0$ as well as the second equality in \eqref{a} we 
obtain that
$$
\Theta_p=-knH\lambda(n,k,H)+k^2\lambda^2(n,k,H).
$$
Then since \eqref{lambda} is equivalent to
$$
-knH\lambda(n,k,H)+k^2\lambda^2(n,k,H)=k(n-k)
$$
we have that equality holds in $(\#)$, and this completes 
the proof of part $(ii)$.

\noindent $(iii)$ Finally, assume that let $H(x)=0$, $k=2$ 
and that $(\#)$ is satisfied for a certain basis. 
From part $(i)$ we have that $p=2$. Moreover, it 
follows from \eqref{ricb}, \eqref{r} and \eqref{ii} that
$\Ric_M(X)=b(n,2,0)=n/2$
for any unit vector $X\in\spa\{e_1,e_2\}$ and
$$
\a_f(X,Y)=\<X,Y\>\eta
$$
if $X,Y\in\spa\{e_1,e_2\}$ and $\eta=\sum_\a\rho_\a\xi_\a$. 
By \eqref{ric} the latter is equivalent to
\be\label{first}
\sum_\a\|A_\a e_i\|^2=(n-2)/2,\;\;i=1,2,
\ee
and therefore, we have 
\be\label{second}
\sum_{j=3}^n\sum_\a\<A_{\a }e_i,e_j\>^2
=(n-2)/2-\|\eta\|^2,\;\;i=1,2.
\ee
It follows that $\|\eta\|\leq\lambda(n,2,0)$, thus 
completing the proof of the direct statement in part $(iii)$.

Conversely, let $\{e_j\}_{1\leq j\leq n}$ be an 
orthonormal basis  of $T_xM$ such that 
$$
\a_f(X,Y)=\<X,Y\>\eta
$$
for any $X,Y\in \spa\{e_1,e_2\}$ where
$\|\eta\|\leq\lambda(n,2,0)$. Moreover, assume that 
$\Ric_M(X)=b(n,2,0)$ for any unit vector 
$X\in\spa\{e_1,e_2\}$. By \eqref{ric} 
we have \eqref{first} and hence \eqref{second}.
Using $H(x)=0$, the second equality 
in \eqref{a} and \eqref{second}  we obtain
$$
\sum_{i=1}^2\sum_{j=3}^n\big(2\|\a_{ij}\|^2-
\<\a_{ii},\a_{jj}\>\big)
=2\sum_{i=1}^2\sum_{j=3}^n\sum_\a\<A_{\a }e_i,e_j\>^2
+4\|\eta\|^2=2(n-2),
$$
which completes the proof of part $(iii)$.\qed

\begin{remark} {\em  
Proposition 5 establishes that the pinching condition ($*$) leads to the 
validity of the Lawson-Simons condition ($\#$). The converse is not true. 
In a recent preprint by the second author \cite{vl}, it was 
proved that a certain pinching condition on the sectional curvature 
implies the Lawson-Simons condition ($\#$). However, there exist 
submanifolds that satisfy the pinching condition on the sectional 
curvature but fail to meet our pinching condition ($*$) on the Ricci 
curvature. Such an example is the $4$-dimensional Veronese 
submanifold (cf. \cite{Itoh}). This is a minimal isometric embedding 
$\psi\colon\R P^4_K\to\Sf^{4+m}$ of the real projective 
space $\R P^4_K$ of sectional curvature $K=2/5$ 
into a unit sphere. This shows that our pinching condition 
($*$) is not equivalent to the Lawson-Simons condition ($\#$). 
}\end{remark}

\section{Proof of Theorem \ref{thm1}}

We first argue that $M^n$ is simply connected, thus orientable, 
by means of an argument given in \cite[p.\ 442]{LS}.
For $p=1$ we have from part $(i)$ of Proposition \ref{prop} 
that the inequality in $(\#)$ is strict at any point 
$x\in M^n$ and for any orthonormal basis of $T_x M$.  
Then Theorem \ref{ls} yields that there are no stable $1$-currents. 
Since in each nontrivial free homotopy 
class there is a length minimizing curve we conclude that 
$\pi_1(M^n)=0$.

According to part $(i)$ of Proposition \ref{prop} the inequality 
$(\#)$ is satisfied for any $1\leq p\leq k$ at any point 
of $M^n$ and for any orthonormal basis of the tangent space at
that point. 
\vspace{1ex}

\noindent $(i)$
At first suppose that 
$H_i(M^n;\Z)=0=H_{n-i}(M^n;\Z)$ for all $1\leq i\leq k$. 
By part $(i)$ of Proposition \ref{prop} and Theorem \ref{ls} 
this is necessarily the case if at some point of $M^n$ the 
inequality $(*)$ is strict. Assume that $k<n/2$. The finitely 
generated abelian group $H_{n-k-1}(M^n;\Z)$ decomposes as
$$
H_{n-k-1}(M^n;\Z)=\Z^{\beta_{n-k-1}(M)}
\oplus {\rm {Tor}}(H_{n-k-1}(M^n;\Z))
$$
where  $\beta_{n-k-1}(M)$ is the $(n-k-1)$-th Betti number of 
$M^n$ and ${\rm {Tor}}(\,)$ is the torsion subgroup. By the 
universal coefficient theorem for cohomology  
(cf. Corollary 4 in \cite[p.\ 244]{Sp}) and Poincar\'e 
duality, we have  
$$
{\rm {Tor}}(H_{n-k-1}(M^n;\Z))\cong{\rm {Tor}}(H^{n-k}(M^n;\Z))
\cong{\rm {Tor}}(H_k(M^n;\Z))=0.
$$ 
Hence $H_{n-k-1}(M^n;\Z)$ is torsion free.
\vspace{1ex}

\noindent $(ii)$
Suppose that $H_i(M^n;\Z)=0=H_{n-i}(M^n;\Z)$ for all 
$1\leq i\leq k$ does not hold. Then  consider the set
$$
S=\left\{1\leq r\leq k\colon H_i(M^n;\Z)
=H_{n-i}(M^n;\Z)=0\;\,\text{for all}\;\,1\leq i\leq r\right\}
$$
which is nonempty because $M^n$ is simply connected and
hence $1\in S$. Set $q=1+\max S$. Then
$H_q(M^n;\Z)\neq 0\;\;\text{or}\;\; H_{n-q}(M^n;\Z)\neq 0$.
By Theorem~\ref{ls}  there exists an orthonormal basis 
$\{e_1,\dots,e_q,e_{q+1},\dots,e_n\}$ of $T_xM$ at any 
point $x\in M^n$ such that equality holds in $(\#)$ for $p=q$.  
Then part $(ii)$ of Proposition \ref{prop} yields that $q=k$ 
and thus 
\be\label{Hk0}
H_k(M^n;\Z)\neq 0 \;\;\text{or}\;\; H_{n-k}(M^n;\Z)\neq0.
\ee
Since $H_{k-1}(M^n;\Z)=0$, 
the universal coefficient theorem for cohomology gives
${\rm {Tor}}(H^k(M^n;\Z))=0$ and Poincar\'e duality that
${\rm {Tor}}(H_{n-k}(M^n;\Z))=0$. Hence,  
$H_{n-k}(M^n;\Z)$ is torsion free and \eqref{Hk0} 
yields $H_k(M^n;\Z)\neq 0$.

If in case $k=2$ we have that $H\neq 0$ everywhere, then part 
$(ii)$ of Proposition \ref{prop} gives that the submanifold 
carries at any point $x$ a Dupin principal normal $\eta$ of multiplicity 
$\ell=\ell(x)\geq k$. Moreover, $\|\eta\|=\lambda(n,k,H)$ and 
the principal normal is collinear with the mean curvature vector 
at  points where it does not vanish. Furthermore, since 
$\|\eta\|=\lambda(n,k,H)$  it follows from \eqref{ric} that
$$
\Ric_M(X)=n-1+nH\lambda(n,k,H)-\lambda^2(n,k,H)
$$ 
for any unit tangent vector $X\in E_\eta$. Since $\lambda(n,k,H)$ 
is a root of the polynomial $P(t)$ given by \eqref{poli}, we obtain 
that $\Ric_M(X)=b(n,k,H)$ for any unit tangent vector $X\in E_\eta$. 

Since $\ell\geq k$, we may assume that $\ell>k$ since otherwise 
we are done. We claim first that $\ell\leq n-k$. Suppose to the 
contrary that $\ell>n-k$. Since $k\leq n/2$ then $\ell>n/2$.  
Consider the isometric immersion $g=i\circ f\colon M^n\to\R^{n+m+1}$ 
where $i\colon\Sf^{n+m}\to\R^{n+m+1}$ is the standard inclusion.
Then $\eta+f$ is a principal normal of $g$ with multiplicity 
$\ell>n/2$. It follows from Theorem $1.23$ in \cite{DT} that
$H_r(M^n;\Z)=0$ for $r$ such that $n-\ell<r<\ell$. In particular, 
we have that $H_k(M^n;\Z)=0$, which is a contradiction that proves 
the claim. 

From the claim, we obtain that $\ell=k$ if $k=n/2$. Then we argue 
that $\ell<n-k$ at any point if  $k<n/2$. Suppose to the contrary 
that there exists a point $x\in M^n$ such that $\ell=n-k$. 
Let $\{e_j\}_{1\leq j\leq n}$ be an orthonormal basis of $T_xM$ 
such that $E_\eta(x)=\spa\{e_{k+1},\dots,e_n\}$. Then we have
\begin{align*}
\sum_{i=1}^k&\sum_{j=k+1}^n\big(2\|\a_{ij}\|^2-
\<\a_{ii},\a_{jj}\>\big)
=-\sum_{i=1}^k\sum_{j=k+1}^n\<\a_{ii},\a_{jj}\>\\
&=-(n-k)\sum_{i=1}^k\<\a_{ii},\eta\>
=-(n-k)\<n\mathcal H-\sum_{j=k+1}^n\a_{jj},\eta\>\\
&=-(n-k)\<n\mathcal H-(n-k)\eta,\eta\>.
\end{align*}
Now a straightforward computation using that 
$\mathcal H=(H/\lambda(n,k,H))\eta$ gives
\begin{align*}
\sum_{i=1}^k&\sum_{j=k+1}^n\big(2\|\a_{ij}\|^2
-\<\a_{ii},\a_{jj}\>\big)-k(n-k)\\
&=\frac{(n-k)(n-2k)}{2k^2}\big(n^2H^2
+nH\sqrt{n^2H^2+4k(n-k)}+2kn\big),
\end{align*}
where the right hand side is positive since $k<n/2$.
Thus for an orthonormal basis of $\{e_j\}_{1\leq j\leq n}$ of $T_xM$ 
such that $E_\eta(x)=\spa\{e_{k+1},\dots,e_n\}$ we have
\be\label{positive}
\sum_{i=1}^k\sum_{j=k+1}^n\big(2\|\a_{ij}\|^2
-\<\a_{ii},\a_{jj}\>\big)-k(n-k)>0.
\ee
On the other hand, for any orthonormal basis of $T_xM$ we obtain
from part $(i)$ of Proposition \ref{prop} that the left hand-side of  
\eqref{positive} is negative, and this is a contradiction.\qed

\section{Proof of Theorem \ref{thm2}}

We have from Theorem \ref{thm1} that
$$
H_i(M^n;\Z)=H_{n-i}(M^n;\Z)=0\;\,\text{for all}\;\,
1\leq i\leq n/2-1.
$$

\noindent $(i)$
Suppose that $H_{n/2}(M^n;\Z)=0$ and hence $M^n$ is a 
homology sphere. Theorem \ref{thm1} also yields that 
$M^n$ is simply connected. Thus in this case the Hurewicz 
homomorphisms between the homotopy and homology groups are 
isomorphisms, and hence $M^n$ is a homotopy sphere. Since the 
generalized Poincar\'e conjecture holds due to the work of 
Smale and Freedman then $M^n$ is homeomorphic to $\Sf^n$. 
\vspace{1ex}

\noindent $(ii)$
Now assume that $H_{n/2}(M^n;\Z)\neq 0$. By part $(i)$ of 
Proposition~\ref{prop} we have that the inequality 
$(\#)$ is satisfied for $p=n/2$ at any $x\in M^n$ 
and for any orthonormal basis of $T_xM$. 
Moreover, by Theorem \ref{ls} at any $x\in M^n$ there 
exists an ordered orthonormal basis 
$\{e_1,\dots,e_{n/2},e_{n/2+1},\dots,e_n\}$ of $T_xM$ such 
that equality holds in $(\#)$ for $p=n/2$. 

We claim that $f$ is minimal. To the contrary, suppose
that $H(x)\neq 0$ at some point $x\in M^n$. Then part $(ii)$ 
of Proposition \ref{prop} gives that the second fundamental 
form of $f$ at $x$ satisfies 
\be\label{sff}
\alpha_f(X,Y)=\<X,Y\>\lambda(n,n/2,H(x))\mathcal{H}(x)/H(x)
\ee
for any $X\in\spa\{e_1,\dots,e_{n/2}\}$ and $Y\in T_xM$. 
It is clear that equality also holds in 
$(\#)$ for $p=n/2$, but now for the ordered orthonormal 
basis $\{e_{n/2+1},\dots,e_n,e_1,\dots,e_{n/2}\}$ of $T_xM$. 
Then \eqref{sff} also holds but now for any 
$X\in\spa\{e_{n/2+1},\dots,e_n\}$ and $Y\in T_xM$. Hence 
$\lambda(n,k,H(x))=H(x)$. But since  
$\lambda(n,n/2,H(x))=H(x)+\sqrt{H^2(x)+1}$ this is a 
contradiction that proves the claim. Since $f$ is minimal 
now the proof follows from the Theorem of Ejiri given 
in \cite{E}.\qed

\section{Proof of Theorem \ref{thm3}}

If $f\colon M^n\to\Sf^{n+m}$ carries  a Dupin principal 
normal vector field $\eta\in\Gamma(N_fM)$ then it is well-known 
that the distribution $x\in M^n\mapsto E_\eta(x)$ is smooth 
and integrable with spherical leaves. The latter condition means 
that there exists $\delta\in\Gamma(E_\eta^\perp)$ such that 
$$
(\n_TS)_{E_\eta^\perp} =\<T,S\>\delta\;\;\text{and}\;\; 
(\n_T\delta)_{E_\eta^\perp}=0
$$
for any $T,S\in\Gamma(E_\eta)$. Moreover, the vector field 
$\eta$ is parallel in the normal connection of $f$ along the 
leaves of $E_\eta$ and $f$ maps each leaf into an 
\mbox{$\ell$-dimensional} round sphere. For details 
cf.\ Proposition $1.22$ of \cite{DT}.

\begin{proposition}\label{sb}
Let $f\colon M^n\to\Sf^{n+m}$, $n\geq 5$, be an isometric 
immersion of a compact manifold $M^n$ of odd dimension 
such that \eqref{ineq2} is satisfied at any point of $M^n$. 
For $n=5$ assume further that $H\neq 0$ at any point. Then 
either $M^n$ is homeomorphic to $\Sf^n$ or the 
following holds:\vspace{1ex}

The submanifold $f$ carries a Dupin principal normal 
vector field $\eta$ with multiplicity $k=(n-1)/2$ such that 
$$
\|\eta\|=\frac{1}{n-1}\big(nH+\sqrt{n^2H^2+n^2-1}\big)
$$ 
and if $X\in E_\eta$ is of unit length then
$$
\Ric_M(X)=\frac{n(n-3)}{n-1}\Big(1+\frac{H}{(n-1)}
\big(nH+\sqrt{n^2H^2+n^2-1}\big)\Big).
$$
Moreover, we have that $M^n$ is diffeomorphic to the total
space of a sphere bundle 
$\Sf^k\hookrightarrow\mathsf{E}\xrightarrow{p} L^{k+1}$ 
having a compact connected manifold as base.
\end{proposition}

\proof By Theorem \ref{thm1} there are just two cases 
to be considered. If case $(i)$ holds then $M^n$ 
is a homology sphere. Since $M^n$ is simply connected, 
arguing as in the proof of Theorem \ref{thm2} it 
follows that $M^n$ is homeomorphic to $\Sf^n$. 

Assume now that part $(ii)$ of Theorem \ref{thm1} holds. 
Hence there is a Dupin principal normal $\eta$ 
with multiplicity $k=(n-1)/2$ at each point of $M^n$. 
We claim that $\eta$ is smooth.  For any normal vector
$\xi$ at $x\in M^n$ we have that the associated shape 
operator satisfies $A_\xi X=\<\eta,\xi\>X$ for any 
$X\in E_\eta(x)$. Thus $\<\eta,\xi\>$ is an eigenvalue
of $A_\xi$ and hence $\eta$ is continuous.  Now that $\eta$
is smooth follows from Theorem $1$ in \cite{Re}.

In the sequel, we identify  
$f$ with the isometric immersion 
$j\circ f\colon M^n\to\R^{n+m+1}$ where 
$j\colon\Sf^{n+m}\to\R^{n+m+1}$ is the standard inclusion 
centered at the origin. Then 
$\mathcal L=f_*E_\eta\oplus\spa\{f_*\delta+\eta-f\}$ is a 
vector subbundle of the vector bundle $f^*T\R^{n+m+1}$ where 
the latter is endowed with the induced metric and connection 
denoted by $\tilde\n$. We claim that the vector subbundle 
$\mathcal L$ is parallel along $E_\eta$. It is easily seen 
that
$$
\tilde\n_T(f_*\delta+\eta-f)=-(1+\|\delta\|^2+\|\eta\|^2)f_*T
$$
and 
$$
\tilde\n_Tf_*S=f_*(\n_TS)_{E_\eta}+\<T,S\>(f_*\delta+\eta-f)
$$
for any $T,S\in\Gamma(E_\eta)$, and the claim follows. The
functions $\|\eta\|$, $\|\delta\|$ and  the vector 
subbundle $\mathcal L$ are constant along each integral leaf 
of $E_\eta$.  Then this is also the case for $r\in C^\infty(M)$ 
given by 
$$
r=1/\sqrt{1+\|\delta\|^2+\|\eta\|^2}
$$
and the map $F\colon M^n\to\R^{n+m+1}$ defined by 
$F=f+r^2(f_*\delta+\eta-f)$.

Let $M_0$ be a maximal $k$-dimensional leaf of $E_\eta$ 
and $f_0\colon M_0\to\R^{n+m+1}$ the isometric immersion 
$f_0=f\circ i$ where $i\colon M_0\to M^n$ is the 
inclusion map. Its second fundamental form is given by 
$$
\alpha_{f_0}(T,S)=\<T,S\>(f_*\delta+\eta-f)
$$
for all $T,S\in\mathfrak{X}(M_0)$. Then $f(M_0)$ is a sphere 
(cf. Proposition 1.22 in \cite{DT}) and $f(M_0)\subset f(x_0)+\mathcal L(x_0)$ for a fixed 
arbitrary point $x_0\in M_0$. Thus $f(M_0)=\Sf^k_{r(x_0)}(F(x_0))$ 
is contained in the affine subspace  $f(x_0)+\mathcal L(x_0)$ 
of $\R^{n+m+1}$ with center $F(x_0)$ and radius $r(x_0)$. 
In addition, the restriction 
$f|_{M_0}\colon M_0\to\Sf^k_{r(x_0)}(F(x_0))$ is a local 
isometry. Since $M_0$ is compact and $\Sf^k_{r(x_0)}(F(x_0))$ 
is simply connected, we have that the restriction $f|_{M_0}$ is 
an isometry between $M_0$ and $f(M_0)=\Sf^k_{r(x_0)}(F(x_0))$. 

Let $L^{k+1}$ be the quotient space of maximal leaves of 
$E_\eta$ and denote by $\pi\colon M^n\to L^{k+1}$ the projection. 
Since $L^{k+1}$ is trivially connected and $M^n$ is 
compact, we obtain from Corollary on p.\ 16 in \cite{P} that it is Hausdorff.  
Thus $L^{k+1}$ is a $(k+1)$-dimensional manifold.  Moreover, the 
projection $\pi$ is a submersion with vertical distribution
$\ker\pi_*=E_\eta$.

That $\mathcal L$ remains constant along each integral 
leaf of the distribution $E_\eta$ gives rise to a vector bundle 
$p\colon\mathsf{V}\to L^{k+1}$ of rank $k+1$ with base space 
$L^{k+1}$ whose fiber over $y\in L^{k+1}$ is 
$\mathsf{V}_y=p^{-1}(y)=\mathcal L(\pi^{-1}(y))$, 
where by parallel translation $\mathcal L(\pi^{-1}(y))$ 
is seen as a linear subspace of $\R^{n+m+1}$.

Let $p\colon\mathsf{E}\to L^{k+1}$ be the unit subbundle  
of the vector bundle $p\colon\mathsf{V}\to L^{k+1}$, that is, 
$$
\mathsf{E}=\left\{(y,v)\colon y\in L^{k+1},v\in\mathsf{V}_y,
\|v\|=1\right\}.
$$
Then let $\Phi\colon M^n\to \mathsf{E}$ be the smooth map
defined by 
$$
\Phi(x)
=\big(\pi(x), \frac{1}{r(x)} (f(x)-F(x))\big)
=\big(\pi(x),r(x)(f(x)-f_*(x)\delta(x)-\eta(x))\big).
$$

We now argue that the map $\Phi$ is a diffeomorphism. At first 
observe that $\Phi$ is injective. This follows from the fact 
that both $r$ and $F$ are constant along each integral leaf of 
$E_\eta$ and that the restriction of $f$ on each leaf is injective. 
We now prove that $\Phi$ is surjective. 
Let $(y,v)\in \mathsf{E}$. Then $y=\pi(x_0)$ 
for some $x_0\in M^n$ and $v\in\mathsf{V}_y=\mathcal L(x_0)$ with 
$\|v\|=1$. Let $M_0$ be the unique maximal leaf of $E_\eta$ containing 
$x_0$. Since $\|v\|=1$, we have that $r(x_0)v$ belongs to the sphere 
$f(M_0)$ centered at $F(x_0)$ with radius $r(x_0)$. Thus there is a 
point $x\in M_0$ such that $r(x_0)v=f(x)-F(x_0)$. Then, we have 
$$
(y,v)=\big(\pi(x_0), \frac{1}{r(x_0)}\big(f(x)-F(x_0)\big)\big)
=\big(\pi(x), \frac{1}{r(x)}\big(f(x)-F(x)\big)\big)
=\Phi(x)
$$
since $x$ and $x_0$ belong to the same leaf of $E_\eta$. Hence 
the map $\Phi$ is a homeomorphism. 

To complete the proof it suffices to show that its differential is 
nonsingular. Indeed, if $\Phi_*(x)X=0$ at some $x\in M^n$, then  
$X\in\ker\pi_*(x)\cap\ker G_*(x)$ where $G=(f-F)/r$. Using that 
$\ker\pi_*(x)=E_\eta(x)$, we obtain that $X$ is vertical and thus 
$r_*(x)X=0$ and $F_*(x)X=0$. Hence, $f_*(x)X=0$ and consequently 
$\ker\Phi_*(x)=\{0\}$. That $\Phi$ is a diffeomorphism 
is a direct consequence of the inverse function theorem.
\qed
\vspace{1.5ex}

Next we  recall several facts relating the cohomology 
of the base space $B$ and the total space $\mathsf{E}$ 
of a sphere bundle 
$\Sf^k\xhookrightarrow{}\mathsf{E}\xrightarrow{p} B$
when $B$ is a compact manifold.

Assume that the sphere bundle is orientable as defined in 
\cite[p.\ 114]{BT}.  Then we have from Theorem 13.2 in \cite{B} 
or \cite[p.\ 438]{Hat} that the cohomology rings of 
$\mathsf{E}$ and $B$ are related by the following long exact 
sequence known as the Gysin sequence:
$$
\cdots\to H^i(B;\Z)\xrightarrow{p^*}H^i(\mathsf{E};
\Z)\xrightarrow{\Phi_i}
H^{i-k}(B;\Z)\xrightarrow{\Psi_i}H^{i+1}(B;\Z)\to\cdots
$$
where all maps are module homomorphisms and $\Psi_i(u)=u\cup e$ 
is the cup product of $u\in H^{i-k}(B;\Z)$ with the Euler 
class $e\in H^{k+1}(B;\Z)$ of the sphere bundle.
\vspace{1ex}

When the coefficient group is the field $\Q$ of rational numbers 
we have the following result.

\begin{lemma}\label{Gys}
Let $\Sf^k\xhookrightarrow{}\mathsf{E}\xrightarrow{p} B$ be 
a sphere bundle such that the total and base 
spaces are both compact connected oriented manifolds. For the  
coefficient group $\Q$ the following assertions hold:\vspace{1ex}

\noindent

\noindent $(i)$ If the fundamental class $[B]$ is in the image of 
$p_*\colon H_*(\mathsf{E})\to H_*(B)$, then the 
cohomology groups of $\mathsf{E}$ and the 
product manifold $B\times\Sf^k$ are isomorphic.\vspace{1ex}

\noindent $(ii)$ If the fiber has even dimension then 
the fundamental class $[B]$ is as in~$(i)$.
\end{lemma}

\proof The proof of part $(i)$ follows from Satz 38 in \cite{G} 
whereas the one of part $(ii)$ from Satz 40 in \cite{G}. 
\vspace{1ex}\qed

\noindent\emph{Proof of Theorem \ref{thm3}:} 
We assume that $M^n$ is not homeomorphic to $\Sf^n$. Then 
Proposition \ref{sb} yields that $M^n$ is diffeomorphic to 
the total space of a sphere bundle 
$\Sf^k\hookrightarrow\mathsf{E}\xrightarrow{p} L$,
$k=(n-1)/2$, with a compact connected manifold as base.

We first argue that $L$ is simply connected. Since 
the fiber bundle  projections are fibrations 
(see Corollary 14 in \cite[p.\ 96]{Sp}) 
then $p\colon\mathsf{E}\rightarrow L$ is 
a fibration. From the homotopy sequence of fibrations,  
we obtain the exact sequence
$$
\pi_1(M^n)\cong\pi_1(\mathsf{E})\xrightarrow{p_*}\pi_1(L)
\rightarrow \pi_0(\Sf^k)=0.
$$
Since $M^n$ is simply connected by Theorem \ref{thm1}, so
is $L$. Then from Propositions $11.2$ and $11.5$ in 
\cite{BT} it follows that the sphere bundle is orientable and, 
consequently, we have the exact Gysin sequence
\be\label{gysin}
\!\!\!\cdots\xrightarrow{\Psi_{i-1}}
H^i(L;\Z)\xrightarrow{p^*}H^i(M;\Z)\xrightarrow{\Phi_i}
H^{i-k}(L;\Z)\xrightarrow{\Psi_i}H^{i+1}(L;\Z)
\xrightarrow{p^*}\cdots
\ee
On the other hand, part $(ii)$ of Theorem \ref{thm1} 
gives that 
$$
H_i(M^n;\Z)=H_{n-i}(M^n;\Z)=0\;\,\text{for all}\;\,1\leq i\leq k-1,
$$
$H_k(M^n;\Z)\cong\Z^{\beta_k(M)}\oplus T$ and 
$H_{k+1}(M^n;\Z)=\Z^{\beta_k(M)}$,  where $T$ is the torsion 
of $H_k(M^n;\Z)\neq0$. Then from the universal 
coefficient theorem for cohomology we have that 
\be\label{new}
H^i(M^n;\Z)=H^{n-i}(M^n;\Z)=0\;\,\text{for all}\;\,1\leq i\leq k-1,
\ee
$H^k(M^n;\Z)=\Z^{\beta_k(M)}$ and 
$H^{k+1}(M^n;\Z)\cong\Z^{\beta_k(M)}\oplus T$. 
It follows from \eqref{gysin} using \eqref{new} that
\be\label{L}
H^i(L;\Z)=0\;\,\text{for all}\;\,1\leq i\leq k-1.
\ee

From \eqref{gysin} for $i=k$ and $i=k+1$ we have 
the exact sequences
\be\label{gysink}
0\to H^k(L;\Z)\xrightarrow{p^*}H^k(M;\Z)\xrightarrow{\Phi_k}
H^0(L;\Z)\xrightarrow{\Psi_k}H^{k+1}(L;\Z)
\xrightarrow{p^*}H^{k+1}(M;\Z)
\ee
and
\be\label{gysink+1}
H^{k+1}(L;\Z)\xrightarrow{p^*}H^{k+1}(M;\Z)
\xrightarrow{\Phi_{k+1}}0.
\ee
In particular, we obtain from \eqref{gysink} that
$p^*\colon H^k(L;\Z)\to H^k(M^n;\Z)$ is a monomorphism. 

Since $L$ is simply connected then
$H^0(L;\Z)\cong\Z\cong H^{k+1}(L;\Z)$. Therefore, 
up to isomorphisms $\Psi_k$ is a homomorphism of $\Z$. Thus either 
$\Psi_k=0$ or up to isomorphisms $\Psi_k$ is a multiplication by 
an integer $q=\Psi_k(1)$.  Therefore, we need to distinguish two cases.
\vspace{1ex}

\noindent\emph{Case $\Psi_k=0$.} We have from \eqref{gysink} that 
$p^*\colon H^{k+1}(L;\Z)\to H^{k+1}(M^n;\Z)$ 
is a monomorphism whereas from \eqref{gysink+1} an 
epimorphism. Hence $H^{k+1}(M^n;\Z)$ and  $H^{k+1}(L;\Z)$ 
are isomorphic and thus $H^{k+1}(M^n;\Z)\cong\Z$. Consequently, 
$H_k(M^n;\Z)\cong\Z\cong H_{k+1}(M^n;\Z)$ and therefore 
the integral homology groups of $M^n$ and $\Sf^k\times\Sf^{k+1}$
are isomorphic.

The exactness of \eqref{gysink} at $H^0(L;\Z)$ gives that 
\be\label{isom}
H^k(M^n;\Z)/{\rm {Im}}\,p^*\cong H^0(L;\Z)\cong\Z.
\ee
Since  ${\rm {Im}}\, p^*\cong p\Z$ for some integer $p$ then
from \eqref{isom} it follows that ${\rm {Im}}\, p^*=0$ and 
since we have seen that $p^*$ is a monomorphism then 
$H^k(L;\Z)=0$. 
This together with \eqref{L} gives that the cohomology ring of 
$L$ is isomorphic to the one of $\Sf^{k+1}$. Thus $L$ 
is a homology sphere. Since $L$  is also simply connected, 
we obtain arguing as in the proof of Theorem \ref{thm2} that 
$L$ is homeomorphic to $\Sf^{k+1}$. Hence $M^n$ 
is diffeomorphic to the total space of a sphere bundle with 
fiber $\Sf^k$ with base homeomorphic to $\Sf^{k+1}$.
\vspace{1ex}

\noindent\emph{Case $\Psi_k\neq 0$.} Then $\Psi_k$ is a 
monomorphism and hence the exactness of \eqref{gysink} at 
$H^0(L;\Z)$ yields that $\Phi_k=0$. Thus 
$p^*\colon H^k(L;\Z)\to H^k(M^n;\Z)$ is an 
epimorphism, and being also a monomorphism, then 
$H^k(L;\Z)\cong H^k(M^n;\Z)\cong\Z^{\beta_k(M)}$. 

The exactness of \eqref{gysink} at $H^{k+1}(L;\Z)$ 
gives that $\ker p^*={\rm {Im}}\,\Psi_k\cong q\Z$ for 
some integer $q\neq 0$ and the exactness of \eqref{gysink+1} 
that ${\rm {Im}}\, p^*=H^{k+1}(M^n;\Z)$. 
It follows that 
$$
H^{k+1}(L;\Z)/\ker\,p^* 
\cong H^{k+1}(M^n;\Z)\cong\Z^{\beta_k(M)}\oplus T 
$$
and thus $\beta_k(M)=0$. Then $H^k(M;\Z)=0$ and 
$H^{k+1}(M^n;\Z)\cong T\cong\Z_q$. Hence the homology group of 
$M^n$ is $\Z_q$ in dimension $k$, $\Z$ in dimensions $0,n$ and 
trivial otherwise. Since $M^n$ is not homeomorphic to $\Sf^n$, 
then \mbox{$q>1$}. Moreover, since 
$p^*\colon H^k(L;\Z)\to H^k(M^n;\Z)$ is a monomorphism then
$H^k(L;\Z)\cong H^k(M^n;\Z)=0$, and we obtain from 
\eqref{L} that $L$ is a homology sphere. As before, it 
follows that $L$ is homeomorphic to $\Sf^{k+1}$ and thus
$M^n$ is diffeomorphic to the total space of a sphere bundle 
with fiber $\Sf^{(n-1)/2}$ with base homeomorphic to 
$\Sf^{(n+1)/2}$.

Finally, let $n=4r+1$. We claim that the homology
of $M^n$ over the integers is isomorphic to that of  
$\Sf^k\times\Sf^{k+1}$. To the contrary assume that this is 
not the case. As shown before, the homology of $M^n$ over 
the integers has to be $\Z_q, q>1$, in dimension $k$, $\Z$ in 
dimensions $0,n$ and trivial otherwise. Then $M^n$ is a 
rational homology sphere, which is in contradiction to
Lemma~\ref{Gys}, and the claim has been proved.

Finally, we recall that by the bundle classification theorem 
given in \S $18.5$ in \cite{St} the equivalence classes of the 
$\Sf^k$ bundles over $\Sf^{k+1}$ are in one-to-one correspondence 
with the homotopy group $\pi_k(SO(k+1))$. If $n=5$ then $M^5$ is 
homeomorphic to $\Sf^2\times\Sf^3$ since any sphere bundle over 
$\Sf^3$ is a product bundle as shown in $\S 26.4$ of \cite{St}. 
If $n=13$, since $\pi_6(SO(7))=0$ (cf. Table 6.VII in 
Appendix A of \cite{MSJ}), the only $6$-sphere bundle over $\Sf^7$ is 
the trivial one. Hence $M^{13}$ is homeomorphic to 
$\Sf^6\times\Sf^7$.\qed

\section{Examples}

The aim of this section is to show that there are many 
 examples of submanifolds, in particular minimal ones, 
for which the inequality $(*)$ is satisfied with equality 
or strictly. 

\subsection{\!\!\!Isoparametric hypersurfaces and  
focal submanifolds}

The subject of isoparametric hypersurfaces in spheres  
introduced by Cartan was continued by Nomizu \cite{N} and then by
M\"{u}nzner \cite{M1,M2}. M\"{u}nzner initiated the geometric 
study of the focal submanifolds of an isoparametric hypersurface 
$f\colon M^n\to\Sf^{n+1}$ with $g$ distinct principal curvatures. 
Any isoparametric hypersurface $M^n$ has associated  just 
two focal submanifolds $M_+$ and $M_-$ whose shape operators
have constant eigenvalues independent of the choices of the 
point and the unit normal vector.

Using that the focal submanifolds of $f$ are obtained as a 
parallel map $f_t$, where $\cot t$ is a principal curvature of 
$f$, M\"{u}nzner computed the shape operator of the focal 
submanifold $f_t$ in terms of the shape operator of $f$ itself. 
In particular, this calculation shows that the assumption that $f$ 
has constant principal curvatures at any point implies that the 
set of eigenvalues for all shape operators are the same with the 
same multiplicity. Then by a symmetry argument, M\"{u}nzner proved 
that if the principal curvatures of $f$ are $\lambda_1>\dots>\lambda_g$,  
with corresponding multiplicities $m_1,\dots,m_g$, then $g$ can 
only be $1,2,3,4$ or $6$ where $m_i=m_{i+2}$ 
(subscripts ${\rm {mod}}\, g$) for
$$
\lambda_i=\cot(\theta+(i-1)\pi/g)\;\;\text{for}\;\;i=1,\dots,g,
$$
and any $\theta\in(0,\pi/2)$. Thus, if $g$ is odd, all of 
the multiplicities must be equal, and if $g$ is even, there are 
at most two distinct multiplicities. Moreover, the one-parameter
family of of isoparametric hypersurfaces contains a unique one
that is minimal.

M\"{u}nzner’s calculation proved further that the focal submanifolds 
must be minimal submanifolds of $\Sf^{n+1}$, as Nomizu \cite{N} had 
shown by a different proof, and that Cartan’s identity is equivalent 
to the minimality of the focal submanifolds. 

If $\cot t$ is not a principal curvature of $f$, then the 
hypersurface $f_t$ is also isoparametric with the $g$ distinct 
principal curvatures $\cot(\theta_i-t)$ for $i=1,\dots,g$. 
If $t=\theta_k\,{\rm {mod}\,}\pi$ then the map $f_t$ is constant 
along each leaf of the $m_k$-dimensional principal foliation $T_k$, 
and the image of $f_t$ is a smooth focal submanifold of $f$ of 
codimension $m_k +1$ in $\Sf^{n+1}$. Then all of the hypersurfaces 
$f_t$ in a family of parallel isoparametric hypersurfaces have 
the same focal submanifolds.

Any compact isoparametric hypersurface must be a level set 
of the restriction to $\Sf^{n+1}\subset\R^{n+2}$ of a homogeneous 
polynomial $F\colon\R^{n+2}\to\R$ of degree $g$ that satisfies the 
Cartan-M\"{u}nzner equations 
$$ 
\|\nabla F\|^2=g^2\|x\|^{2g-2},
\;\;\Delta F=\frac{1}{2}(m_2-m_1)g^2\|x\|^{2g-2}.
$$
These are called the Cartan-M\"{u}nzner polynomials and 
$h=F|_{\Sf^{n+1}}$ takes values in $[-1,1]$. For each $-1<t<1$, 
the level set $h^{-1}(t)$ is an isoparametric hypersurface. 
The level sets $M_+^{n-m_1-1}=h^{-1}(1)$ and 
$M_-^{n-m_2-1}=h^{-1}(-1)$ are the two focal submanifolds 
being both minimal submanifolds in $\Sf^{n+1}$ with non-flat 
normal bundle (cf.\ \cite{QT}). 

If $g=2$ the isoparametric hypersurfaces in $\Sf^{n+1}$ are
just the generalized Clifford torus $\mathbb T_p^n(r)$ and the 
focal submanifolds are $\Sf^p$ and $\Sf^{n-p}$, respectively. 
If $g=3$ Cartan showed that the multiplicities must be 
$m_1=m_2=m=1, 2, 4$ or $8$ and the focal submanifolds are 
Veronese embeddings of $\mathbb FP^2$ in $\Sf^{3m+1}$, 
where $\mathbb F=\R,\C,\mathbb H,\mathbb O$ corresponding 
to $m=1,2,4,8$. Meanwhile, the hypersurfaces are tubes of 
constant radius over the focal submanifolds. The isoparametric 
hypersurfaces with $g=3$ are homogeneous and are called 
Cartan hypersurfaces. 

If $g=4$, the isoparametric hypersurfaces are either 
homogeneous with $(m_1,m_2)=(2,2),(4,5)$ or of FKM-type. 
The ones of FKM-type are defined as follows. Let ${P_0,\dots,P_r}$ 
be asymmetric Clifford system on $\R^{2\ell}$, that is, the  
$P_i$'s are symmetric matrices satisfying 
$P_iP_j +P_jP_i=2\delta_{ij}I_{2\ell}$. Ferus, 
Karcher and M\"{u}nzner \cite{FKM} considered the 
Cartan-M\"{u}nzner polynomials $F\colon\R^{2\ell}\to\R$ 
of degree $4$ given by 
$$
F(x)=\|x\|^4-2\sum_{i=0}^r\<P_ix,x\>^2.
$$
Then the level set of the function $F|_{\Sf^{2\ell-1}}$ 
is an isoparametric hypersurface in $\Sf^{2\ell-1}$ called 
of \emph{FKM-type}.
The multiplicity pair is $(m_1,m_2)=(r,s\delta_r-r-1)$ 
provided $r,s\in\mathbb{N}$ and $s\delta_r-r-1 > 0$, where 
$\delta_r$ is the dimension of the irreducible module of 
the Clifford algebra $C_{r-1}$. For a table with the possible 
values of $(m_1,m_2)$ we refer to \cite{C}.

If $g=6$ then the multiplicities satisfy $m_1=m_2=m=1$ or 
$2$ and the hypersurfaces are homogeneous. 
\vspace{1ex}

Next we show that the Cartan minimal hypersurfaces 
of dimensions $12$ and $24$ satisfy the inequality $(*)$ 
for appropriate values of $k$.

\begin{proposition}\label{isop3}
The following facts hold:\vspace{1ex}

\noindent $(i)$ The Cartan minimal hypersurface 
$f\colon M^{12}\to\Sf^{13}$ satisfies  
$(*)$ for at most $k=3$ and in that case with equality. 
Thus the second fundamental form has the structure given 
by part $(ii)$ of Theorem \ref{thm1} 
with $\lambda (12,3,0)=\sqrt{3}$ a principal curvature 
with multiplicity $4$. 
\vspace{2ex}

\noindent $(ii)$ The Cartan minimal hypersurface 
$f\colon M^{24}\to\Sf^{25}$ satisfies  
$(*)$ for at most $k=6$ and in that case with equality. 
Thus the second fundamental form has the structure given 
by part $(ii)$ of Theorem \ref{thm1} with 
$\lambda (24,6,0)=\sqrt{3}$ a principal curvature with 
multiplicity $8$. 
\end{proposition}

\proof An isoparametric hypersurface with $g=3$ has
the principal curvatures 
\begin{align*}
\lambda_1&=\cot\theta,\:\:0<\theta<\pi/2,\\
\lambda_2&=\cot(\theta+\pi/3)=\frac{\cot\theta
-\sqrt{3}}{\sqrt{3}\cot\theta+1}
=\frac{\lambda_1-\sqrt{3}}{\sqrt{3}\lambda_1+1},\\
\lambda_3&=\cot(\theta+2\pi/3)
=\frac{\cot\theta+\sqrt{3}}{1-\sqrt{3}\cot\theta}
=\frac{\lambda_1+\sqrt{3}}{1-\sqrt{3}\lambda_1}
\end{align*}
all with multiplicity $m$. Thus a Cartan hypersurface 
hypersurface is minimal if $\lambda_1+\lambda_2+\lambda_3=0$. 
Hence it is minimal if the principal curvatures are 
$\sqrt{3},0,-\sqrt{3}$ with multiplicities $4$ for $n=12$ 
and $8$ for $n=24$. 
Then the maximum $k$ for which $(*)$ is satisfied is $k=3$ 
if $n=12$ and and $k=6$ if $n=24$, with equality in both 
cases. Moreover, we have that 
$\lambda (12,3,0)=\sqrt{3}=\lambda (24,6,0)$. 
That the second fundamental form  has the structure given 
by part $(ii)$ of Theorem \ref{thm1} follows from the fact 
that $\sqrt{3}$ is a principal curvature with multiplicity 
four when $n=12$ and eight when $n=24$.
\vspace{1ex}\qed

Next we provide several examples of minimal 
isoparametric hypersurfaces with $g=4$ that satisfy 
the inequality $(*)$ for appropriate values of $k$.  

\begin{proposition}\label{isop} All the minimal isoparametric 
hypersurfaces with $g=4$ and multiplicity pair $(m_1,m_2)$ 
given next satisfy $(*)$ with strict inequality and 
with the exception of case $(i)$ are of FKM-type.
\vspace{1ex}

\noindent $(i)$ The homogeneous hypersurface
$f\colon M^{18}\to\Sf^{19}$ with $(m_1,m_2)=(4,5)$ 
satisfies $(*)$ only for $k=2$.
\vspace{2ex}

\noindent $(ii)$ The hypersurface 
$f\colon M^{14}\to\Sf^{15}$  with 
$(m_1,m_2)=(4,3)$ satisfies $(*)$ only for $k=2$.
\vspace{2ex}

\noindent $(iii)$ The hypersurface 
$f\colon M^{22}\to\Sf^{23}$  with 
$(m_1,m_2)=(4,7)$ satisfies $(*)$ only for $k=2$. 
\vspace{2ex}

\noindent $(iv)$ The hypersurface 
$f\colon M^{30}\to\Sf^{31}$  with $(m_1,m_2)=(4,11)$ 
satisfies $(*)$ at most for $k=3$.
\vspace{2ex}

\noindent $(v)$ The hypersurface 
$f\colon M^{30}\to\Sf^{31}$  with  $(m_1,m_2)=(6,9)$ 
satisfies $(*)$ at most for $k=3$.
\vspace{2ex}

\noindent $(vi)$ The hypersurface 
$f\colon M^{8m-4}\to\Sf^{8m-3}$ with 
$(m_1,m_2)=(4,4m-5)$ and $m\geq 2$ satisfies $(*)$ 
only for $k=2$.
\vspace{1ex}

\noindent $(vii)$ The hypersurface 
$f\colon M^{2s\delta_r-2}\to\Sf^{2s\delta_r-1}$ of FKM 
type with  $(m_1,m_2)=(r,s\delta_r-r-1)$ satisfies  
$(*)$ for $s$ large enough and at most for $k=r/2$.
\end{proposition}

\proof The principal curvatures are given by
\begin{align*}
\lambda_1&=\cot\theta,\:\:0<\theta<\pi/2,\;\;
\lambda_2=\cot(\theta+\pi/4)
=\frac{\cot\theta-1}{\cot\theta+1}
=\frac{\lambda_1-1}{\lambda_1+1},\\
\lambda_3
&=\cot(\theta+\pi/2)=-1/\lambda_1\;\;\text{and}\;\;
\lambda_4=\cot(\theta+3\pi/4)=-1/\lambda_2
\end{align*}
with corresponding multiplicities $m_1,m_2,m_1,m_2$. 
Hence, the hypersurface is minimal if 
$m_1(\lambda_1+\lambda_3)+m_2(\lambda_2+\lambda_4)=0$
which holds if and only if
$$
\lambda_1=\sqrt{m_2/m_1}
+\sqrt{1+m_2/m_1}.
$$
We have from Proposition $3.1$ in \cite{QTY} and its 
proof that the Ricci curvature is positive if and only 
if $m_1,m_2\geq 2$ and satisfies that
$\min{\rm {Ric}}=2(m_1+m_2)-1-\max\{\lambda_1^2,\lambda_4^2\}$.
Thus $(*)$ takes the form
$$
2(m_1+m_2)-1-\max\{\lambda_1^2,\lambda_4^2\}\geq n(k-1)/k.
$$

\noindent $(i)$ The homogeneous hypersurface 
$f\colon M^{18}\to\Sf^{19}$ with $(m_1,m_2)=(4,5)$ has 
principal curvatures 
$$
\lambda_1=(3+\sqrt{5})/2\;\;\text{and}\;\;\lambda_4=-\sqrt{5}.
$$ 
Then the Ricci curvature satisfies $(*)$ only for $k=2$
and with strict inequality.
\vspace{2ex}

\noindent $(ii)$ The  hypersurface 
$f\colon M^{14}\to\Sf^{15}$  with  
$(m_1,m_2)=(4,3)$ has principal curvatures 
$$
\lambda_1=(\sqrt{3}+\sqrt{7})/2,\;\;\text{and}\;\;
\lambda_4=-(\sqrt{3}+\sqrt{7}+2)/(\sqrt{3}+\sqrt{7}-2).
$$ 
Then  Ricci curvature satisfies $(*)$ only for 
$k=2$ and with strict inequality.  
\vspace{2ex}

\noindent $(iii)$ The  hypersurface 
$f\colon M^{22}\to\Sf^{23}$  with  
$(m_1,m_2)=(4,7)$  has principal curvatures 
$$
\lambda_1=(\sqrt{7}+\sqrt{11})/2\;\;\text{and}\;\; 
\lambda_4=-(\sqrt{7}+\sqrt{11}+2)/(\sqrt{7}+\sqrt{11}-2). 
$$
The Ricci curvature satisfies $(*)$ only 
for $k=2$ and with strict inequality. 
\vspace{2ex}

\noindent $(iv)$ The hypersurface 
$f\colon M^{30}\to\Sf^{31}$  with 
$(m_1,m_2)=(4,11)$ has principal curvatures 
$$
\lambda_1=(\sqrt{11}+\sqrt{15})/2\;\;\text{and}\;\; 
\lambda_4=-(\sqrt{11}+\sqrt{15}+2)/(\sqrt{11}+\sqrt{15}-2). 
$$
The Ricci curvature satisfies $(*)$ only for $k=2,3$ 
and with strict inequality.
\vspace{2ex}

\noindent $(v)$ The  hypersurface 
$f\colon M^{30}\to\Sf^{31}$  with  
$(m_1,m_2)=(6,9)$  has principal curvatures 
$$
\lambda_1=(\sqrt{3}+\sqrt{5})/\sqrt{2}\;\;\text{and}\;\; 
\lambda_4
=-(\sqrt{3}+\sqrt{5}+\sqrt{2})/(\sqrt{3}+\sqrt{5}-\sqrt{2}). 
$$
The Ricci curvature satisfies $(*)$ only for $k=2,3$ and 
with strict inequality.
\vspace{2ex}

\noindent $(vi)$ The  hypersurface $f\colon M^{8m-2}\to\Sf^{8m-1}$ 
with $(m_1,m_2)=(4,4m-5)$ and $m\geq 2$ has principal curvatures 
\begin{align*}
\lambda_1&=(\sqrt{4m-5}+\sqrt{4m-1})/2,\\
\lambda_4&=-(\sqrt{4m-5}+\sqrt{4m-1}+2)/(\sqrt{4m-5}+\sqrt{4m-1}-2). 
\end{align*}
The Ricci curvature satisfies $(*)$ only for $k=2$ and
with strict inequality.
\vspace{1ex}

\noindent $(vii)$ The hypersurface 
$f\colon M^{2s\delta_r-2}\to \Sf^{2s\delta_r-1}$ with 
$(m_1,m_2)=(r,s\delta_r-r-1)$ has principal 
curvatures 
\begin{align*}
\lambda_1&=(\sqrt{s\delta_r-r-1}+\sqrt{s\delta_r-1})/\sqrt{r},\\
\lambda_4&=-(\sqrt{s\delta_r-r-1}+\sqrt{s\delta_r-1}+\sqrt{r})
/(\sqrt{s\delta_r-r-1}+\sqrt{s\delta_r-1}-\sqrt{r}) 
\end{align*}
and hence $\lambda_1^2>\lambda_4^2$ holds for $s$ large enough. 
We have that $(*)$ has the form
$k\leq\ 2(s\delta_r-1)/(1+\lambda_1^2)$.
Then for $s$ large enough we obtain that the Ricci curvature
satisfies $(*)$ with strict inequality for $2\leq k\leq r/2$.
\qed

\begin{remark} {\em  From a result in \cite{PT} we have that 
there are Dupin hypersurfaces that satisfy $(*)$ which are not 
isoparametric. They are obtained by deforming the FKM-type 
isoparametric hypersurfaces given by Proposition~\ref{isop}.
}\end{remark}

The following results show that the focal submanifolds of 
infinitely many isoparametric hypersurfaces with $g=4$ do 
satisfy  the inequality $(*)$ for appropriate values of $k$.

\begin{proposition}\label{focal} All of the following focal 
submanifolds $M_+$ of isoparametric hypersurfaces in spheres 
with $g=4$ are minimal and satisfy $(*)$ with strict inequality.
\vspace{1ex}

\noindent $(i)$ If $(m_1,m_2)=(1,m-2)$, $m\geq 5$, then  
$\psi\colon M_+^{2m-3}\to\Sf^{2m-1}$ satisfies $(*)$ for 
at most $k=(2m-3)/3=\dim M_+/3$ and 
$$
H_i(M_+^{2m-3};\Z)=
\begin{cases}
\,\Z&\text{if\, } i=0,m-2,m-1,2m-3 \\[1mm]
\,0&\text{if otherwise.}
\end{cases}
$$

\noindent $(ii)$ If $(m_1,m_2)=(2,2m-1)$, $m\geq 2$, then  
$\psi\colon M_+^{4m}\to\Sf^{4m+3}$ satisfies $(*)$
for at most $k=m=\dim M_+/4$ and
$$
H_i(M_+^{4m};\Z)=
\begin{cases}
\,\Z&\text{if\, } i=0,2m-1,2m+1,4m \\[1mm]
\,0&\text{if otherwise.}
\end{cases}
$$

\noindent $(iii)$ If $(m_1,m_2)=(4,4m-5), m\geq 3$, then 
$\psi\colon M_+^{8m-6}\to\Sf^{8m-1}$ satisfies $(*)$ 
for at most $k\leq (4m-3)/3=\dim M_+/6$ 
and 
$$
H_i(M_+^{8m-6};\Z)=
\begin{cases}
\,\Z&\text{if\, } i=0,4m-5,4m-1,8m-6 \\[1mm]
\,0&\text{if otherwise.}
\end{cases}
$$

\noindent $(iv)$  If $(m_1, m_2)=(6,9)$, then  
$\psi\colon M_+^{24}\to\Sf^{31}$ satisfies $(*)$ 
only for at most $k=3$ and
$$
H_i(M_+^{24};\Z)=
\begin{cases}
\,\Z&\text{if\, } i=0,9,15,24 \\[1mm]
\,0&\text{if otherwise.}
\end{cases}
$$

\noindent $(v)$  If $(m_1,m_2)=(4,5)$, then 
$\psi\colon M_+^{14}\to\Sf^{19}$ satisfies $(*)$ 
with strict inequality only if $k=2$ and 
$$
H_i(M_+^{14};\Z)=
\begin{cases}
\,\Z&\text{if\, } i=0,5,9,14 \\[1mm]
\,0&\text{if otherwise.}
\end{cases}
$$

\noindent $(vi)$  If the isoparametric hypersurface is of 
FKM-type with multiplicity pair $(m_1,m_2)=(r,s\delta_r-r-1)$,  
then $\psi\colon M_+^{2s\delta_r-r-2}\to\Sf^{2s\delta_r-1}$ 
satisfies $(*)$ with strict inequality only for at most
$k=(2s\delta_r-r-2)/(r+2)=\dim M_+/(r+2)$ and 
$$
H_i(M_+^{2s\delta_r-r-2};\Z)=
\begin{cases}
\,\Z&\text{if\, } i=0,s\delta_r-r-1,s\delta_r,2s\delta_r-r-2 \\[1mm]
\,0&\text{if otherwise.}
\end{cases}
$$
\end{proposition}

\proof According to part $(i)$ of Section 4.5 in \cite{FKM} the 
principal curvatures of the focal submanifold 
$\psi\colon M_+^{m_1+2m_2}\to\Sf^{2(m_1+m_2)+1}$ are $0$ with 
multiplicity $m_1$ and $1,-1$ with multiplicity $m_2$. 
We have from \eqref{ric} that the Ricci curvature of the focal 
submanifold is $\Ric>2(m_2-1)$. Then $(*)$ is satisfied in 
cases $(i)$ to $(vi)$  if and only if 
$2\leq k\leq(m_1+2m_2)/(m_1+2)$
and always with strict inequality. According to Theorem \ref{thm1}, 
all the above focal submanifolds are oriented, and thus their 
cohomology groups over the integers according to Satz $5$ 
in \cite{M2} are given by
$$
H^i(M_+^{m_1+2m_2};\Z)=
\begin{cases}
\,\Z&\text{if\, } i=0,m_2,m_1+m_2,m_1+2m_2 \\[1mm]
\,0&\text{if otherwise.}
\end{cases}
$$
In particular, the cohomology groups are torsion free. By the 
universal coefficient theorem for cohomology we have that the 
torsion subgroups of $H^{i+1}(M;\Z)$ and $H_i(M;\Z)$ are isomorphic 
for any $i$. In addition, the free subgroups of $H_i(M;\Z)$ and
$H^i(M;\Z)$ are isomorphic for any $i$. Thus the homology 
groups of the above focal submanifolds are also torsion free 
and therefore  $H_i(M;\Z)$ is isomorphic to $H^i(M;\Z)$ for 
any $i$. Thus the homology groups are as in 
the statement.\qed

\subsection{Embeddings of projective spaces}

Let $\mathbb FP^m$ denote the projective space over 
$\mathbb F=\R,\C,\Hy$ when considered as the quotient 
space of the unit sphere 
$\Sf^{(m+1)d-1}=\{x\in\mathbb{F}^{m+1}\colon \|x\|=1\}$
obtained by identifying $x$ with $\lambda x$ where 
$|\lambda|=1$ and $d=1,2$ or $4$ if $F=\R,\C$ or $\Hy$, 
respectively. 
Let $\mathbb FP^m$ be endowed with the canonical Riemannian 
metric such that $\pi\colon\Sf^{(m+1)d-1}\to\mathbb{F}P^m$ 
is a Riemannian submersion; see \cite{Sk} for details for 
this and the sequel. Then these Einstein manifolds admit 
isometric minimal embeddings 
$\psi_m\colon\mathbb FP^m\to\Sf^{\frac{1}{2}m(m+1)d+m-1}$.
In addition, we also consider the isometric minimal embedding 
$\psi\colon\mathbb OP^2\to\Sf^{25}$ of the Cayley plane 
$\mathbb OP^2$ equipped with the canonical metric of 
$Q$-sectional curvature $4/3$. 
The homology groups are 
$$
H_i(\mathbb FP^m;\Z)=
\begin{cases}
\,\Z&\text{if\, } i=0 \;{\rm {mod}}\, d \\[1mm]
\,0&\text{if otherwise}, \\
\end{cases}
$$
where $d=1,2,4,8$ if $\mathbb F=\R,\C,\Hy,\mathbb O$, 
respectively.

\begin{proposition}\label{proj} Among the submanifolds 
given above the ones that satisfy the inequality $(*)$ 
are the following:\vspace{1ex}

\noindent $(i)$ The embedding 
$\psi_m\colon\C P^m_{2m/(m+1)}\to\Sf^{m^2+2m-1}$ satisfies $(*)$ 
only for $k=2$ and with equality since in this case
$\Ric=b(2m,2,0)$.
\vspace{1ex}

\noindent $(ii)$ The embedding 
$\psi_m\colon\mathbb HP^m\to\Sf^{2m^2+3m-1}$ satisfies 
$(*)$ at most for $k=3$ if $m=2$ and only for $k=2$ if 
$m>2$.  Moreover, equality holds only in the first 
case since $\Ric=b(4m,k,0)$.
\vspace{1ex}

\noindent $(iii)$ The embedding 
$\psi\colon\mathbb OP^2\to\Sf^{25}$ satisfies $(*)$
at most for $k=4$ and equality holds since 
$\Ric=b(16,k,0)$. 
\end{proposition}

\proof Notice that the minimal embeddings of 
$\R P^m$ cannot satisfy $(*)$.\vspace{1ex}

\noindent $(i)$ The Ricci curvature $\C P^m$ is $m$ and 
the result follows.
\vspace{1ex}

\noindent $(ii)$ We obtain from \cite{ISH} that its scalar 
curvature is $8m^2(m+2)/(m+1)$. Hence its Ricci curvature 
is $2m(m+2)/(m+1)$, and then the result follows easily. 
\vspace{1ex}

\noindent $(iii)$ We have from Remark $1$ in \cite{Tan} 
that the scalar curvature is $192$. Hence the Ricci curvature 
is $12$, and then the result follows easily.\qed


\bigskip

\noindent{\bf Funding} Marcos Dajczer is  partially supported 
by the grant PID2021-124157 NB-I00 funded by 
MCIN/AEI/10.13039/501100011033/ `ERDF A way of making Europe',
Spain, and is also supported by Comunidad Aut\'{o}noma de la Regi\'{o}n
de Murcia, Spain, within the framework of the Regional Programme
in Promotion of the Scientific and Technical Research (Action Plan 2022),
by Fundaci\'{o}n S\'{e}neca, Regional Agency of Science and Technology,
REF, 21899/PI/22.

\medskip
\noindent{\bf Data availability} No datasets were generated or analysed during the current study.
\medskip

\noindent
{\bf Declarations}
\medskip
\noindent
{\bf Conflict of interest} The authors declare no competing interests.

\noindent Marcos Dajczer\\
Departamento de Matemáticas\\ 
Universidad de Murcia, Campus de Espinardo\\ 
E-30100 Espinardo, Murcia, Spain\\
e-mail: marcos@impa.br
\bigskip

\noindent Theodoros Vlachos\\
University of Ioannina \\
Department of Mathematics\\
45110 Ioannina -- Greece\\
e-mail: tvlachos@uoi.gr
\end{document}